\documentclass[final]{sig-alternate}

\usepackage{amsmath}

\usepackage{amsfonts}

\usepackage{amssymb}

\usepackage{stmaryrd}

\usepackage[mathscr]{eucal}

\usepackage{amsbsy}

\usepackage{bm}

\usepackage{bbding}

\usepackage{array}

\usepackage[usenames]{color}

\usepackage{tabularx}

\usepackage{paralist}

\usepackage{fancyvrb}

\usepackage[position=below]{subfig}

\usepackage{boxedminipage}

\usepackage{multirow}
\usepackage{algorithmic}

\usepackage{graphicx}
\usepackage{ifpdf}
\ifpdf
\DeclareGraphicsExtensions{.pdf,.png}
\else
\DeclareGraphicsExtensions{.eps}
\fi

\usepackage[draft,hypertexnames=false]{hyperref}

\usepackage[notref,notcite]{showkeys}

\newcommand{\Sec}[1]{\hyperref[sec:#1]{\S\ref*{sec:#1}}} 
\newcommand{\Section}[1]{\hyperref[sec:#1]{Section~\ref*{sec:#1}}} 
\newcommand{\App}[1]{\hyperref[sec:#1]{\ref*{sec:#1}}} 
\newcommand{\Eqn}[1]{\hyperref[eq:#1]{(\ref*{eq:#1})}} 
\newcommand{\Fig}[1]{\hyperref[fig:#1]{Figure~\ref*{fig:#1}}} 
\newcommand{\Tab}[1]{\hyperref[tab:#1]{Table~\ref*{tab:#1}}} 
\newcommand{\Thm}[1]{\hyperref[thm:#1]{Theorem~\ref*{thm:#1}}} 
\newcommand{\Cor}[1]{\hyperref[cor:#1]{Corollary~\ref*{cor:#1}}} 
\newcommand{\Alg}[1]{\hyperref[alg:#1]{Algorithm~\ref*{alg:#1}}} 
\newcommand{\Def}[1]{\hyperref[def:#1]{Definition~\ref*{def:#1}}} 

\newcommand{\Real}{{\mathbb R}}

\newenvironment{inlinemath}{$}{$}

\newcommand{\Tra}{^{{\sf T}}} 
\newcommand{\VC}[1]{{{\sf vec} {\left(#1 \right)}}} 
\newcommand{\V}[1]{{\bm{\mathbf{\MakeLowercase{#1}}}}} 
\newcommand{\Vhat}[1]{{\bm \hat{\mathbf{\MakeLowercase{#1}}}}} 

\newcommand{\M}[1]{{\bm{\mathbf{\MakeUppercase{#1}}}}} 
\newcommand{\Mhat}[1]{{\bm{\hat \mathbf{\MakeUppercase{#1}}}}} 
\newcommand{\MC}[2]{\V{#1}_{#2}} 
\newcommand{\MhatC}[2]{\Vhat{#1}_{#2}} 
\newcommand{\Mn}[2]{\M{#1}^{(#2)}} 
\newcommand{\MnE}[3]{\MakeLowercase{#1}^{(#2)}_{#3}} 
\newcommand{\MCTra}[2]{\V{#1}^{{\sf T}}_{#2}} 

\newcommand{\Kron}{\otimes} 
\newcommand{\Khat}{\odot} 
\newcommand{\Hada}{\ast} 

\newcommand{\T}[1]{\boldsymbol{\mathscr{\MakeUppercase{#1}}}} 
\newcommand{\That}[1]{\boldsymbol{\hat \mathscr{\MakeUppercase{#1}}}} 
\newcommand{\Tbar}[1]{\boldsymbol{\bar \mathscr{\MakeUppercase{#1}}}} 
\newcommand{\TM}[2]{\M{#1}_{(#2)}} 
\newcommand{\TE}[2]{\MakeLowercase{#1}_{#2}} 
\newcommand{\TEP}[2]{\left( #1 \right)_{#2}} 
\newcommand{\Mz}[2]{\TM{#1}{#2}} 

\newcommand{\KOp}[1]{\llbracket #1 \rrbracket} 

\newcommand{\SizeN}[2]{{#1}_1 \times {#1}_2 \times \cdots \times {#1}_{#2}}

\newcommand{\SubscriptN}[2]{{#1}_1 {#1}_2 \cdots {#1}_{#2}}

\newcommand{\SumN}[3]{\sum_{{#1}_1=1}^{{#2}_1} %
  \sum_{{#1}_2=1}^{{#2}_2} %
  \cdots %
  \sum_{{#1}_{#3}=1}^{{#2}_{#3}}}

\newcommand{\norm}[1]{\left\lVert \, #1 \, \right\rVert}

\newcommand{\ip}[2]{\langle \, #1,#2 \, \rangle}

\newcommand{\FD}[2]{\frac{\partial #1}{\partial #2}}

\begin{document}

\title{All-at-once Optimization for\\ Coupled Matrix and Tensor Factorizations}
\numberofauthors{3}
\author{
\alignauthor
Evrim Acar\\
       \affaddr{Faculty of Life Sciences, University of Copenhagen}\\
       \email{evrim@life.ku.dk}
\alignauthor
Tamara G. Kolda\\
       \affaddr{Sandia National Laboratories}\\
       \affaddr{Livermore, CA 94551-9159}\\
       \email{tgkolda@sandia.gov}
\alignauthor
Daniel M. Dunlavy\\
       \affaddr{Sandia National Laboratories}\\
       \affaddr{Albuquerque, NM 87185-1318}\\
       \email{dmdunla@sandia.gov}
}

\makeatletter
\let\@copyrightspace\relax
\makeatother

\maketitle
\begin{abstract}
 Joint analysis of data from multiple sources has the potential to improve
 our understanding of the underlying structures in complex data sets.
 For instance, in restaurant recommendation systems, recommendations can be based on
 rating histories of customers. In addition to rating histories,
 customers' social networks (e.g., Facebook friendships) 
 and restaurant categories information (e.g., Thai or Italian)
 can also be used to make better recommendations. The task of fusing data, however, is challenging
 since data sets can be incomplete and heterogeneous, i.e., data consist of both
 matrices, e.g., the \emph{person} by \emph{person} social network matrix or the \emph{restaurant} by \emph{category} matrix,
 and higher-order tensors, e.g., the ``ratings'' tensor of the form \emph{restaurant} by \emph{meal} by \emph{person}.

 In this paper, we are particularly interested in fusing data sets with the goal of capturing their underlying latent
 structures. We formulate this problem as a coupled matrix and tensor factorization (CMTF) problem where
 heterogeneous data sets are modeled by fitting outer-product models to higher-order tensors and matrices
 in a coupled manner. Unlike traditional approaches solving this problem using alternating algorithms, we propose an all-at-once
 optimization approach called CMTF-OPT (CMTF-OPTimization), which is a gradient-based optimization approach for joint analysis of matrices
 and higher-order tensors. We also extend the algorithm to handle coupled incomplete data sets. Using numerical experiments, we demonstrate that the proposed all-at-once approach is more accurate than the alternating least squares approach.
\end{abstract}



\section{Introduction}
\label{sec:introduction}

With the ability to access massive amounts of data as a result of recent
technological advances, e.g., the Internet, communication and multi-media
devices, genomic technologies and  new medical diagnostic techniques,
we are faced with data sets from multiple sources. For instance, in restaurant
recommendation systems, online review sites like Yelp have access to shopping
histories of customers, friendship networks of those customers, as well as categorizations
of the restaurants. Similarly, for medical diagnoses,
several types of data are collected from a patient; for example, EEG
(electroencephalogram) and ECG (electrocardiogram) monitoring data, fMRI (functional
Magnetic Resonance Imaging) scans, and other data gathered from laboratory tests.

\begin{figure}[t]
\centering
\includegraphics[width=2.5in]{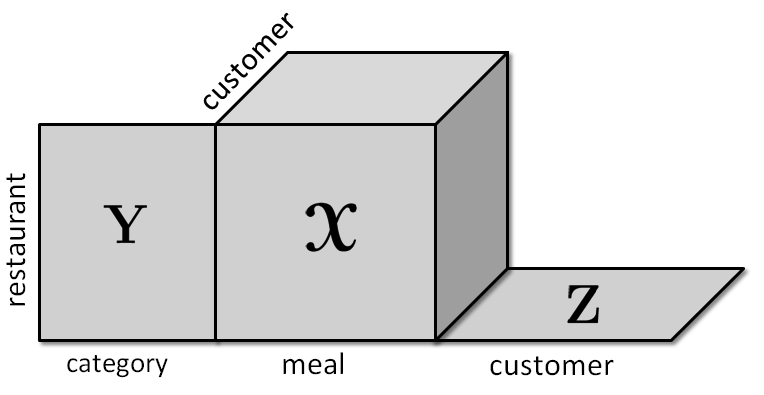}
\caption[Coupled Data Sets of Different Orders]{%
Coupled Data Sets of Different Orders. Each tensor entry indicates the rating of a customer for a specific meal (i.e., breakfast, lunch, dinner) at a particular restaurant. Matrices
$\M{Y}$ and $\M{Z}$ show restaurant categories (i.e., Thai, Chinese, Italian) and social network information, respectively.
Thus, the third-order tensor $\T{X}$ of type \emph{restaurant} by \emph{meal} by \emph{customer}
can be coupled with the \emph{restaurant} by \emph{category}  matrix $\M{Y}$ and the \emph{customer} by \emph{customer} matrix $\M{Z}$.
} \label{fig:CoupledEx}
\end{figure}

Analysis of data from multiple sources requires handling of data sets
of different orders \cite{BaBaMe07,MeSm10,SmWeBo00},
i.e., matrices and/or higher-order tensors. For instance, Banerjee  et~al.~\cite{BaBaMe07},
discusses the problem of analyzing heterogeneous data sets with a goal of simultaneously
clustering different classes of entities based on multiple relations, where
each relation is represented as a matrix (e.g., \emph{movies} by \emph{review words} matrix showing
movie reviews) or a higher-order tensor (e.g., \emph{movies} by \emph{viewers}
by \emph{actors} tensor showing viewers' ratings). Similarly, coupled analysis of matrices
and tensors has been a topic of interest in the areas of community detection \cite{LiSuCaKo09},
collaborative filtering \cite{ZhCaZhXiYa10}, and chemometrics \cite{SmWeBo00}.

As an example of data sets from multiple sources, without loss of generality, suppose
we have a third-order tensor, $\T{X} \in \Real^{I \times J \times K}$,
and a matrix, $\M{Y} \in \Real^{I \times M}$, coupled in the first dimension (mode) of each.
The common latent structure in these data sets can be extracted through coupled
matrix and tensor factorization (CMTF), where an $R$-component CMTF model of a tensor $\T{X}$ and
a matrix $\M{Y}$ is defined as:
\begin{multline}
\label{eq:fABCV}
   f(\M{A},\M{B},\M{C},\M{V})=  \norm{\T{X} - \KOp{\M{A},\M{B},\M{C}}}^2 + \norm{\M{Y} - \M{A}\M{V}\Tra}^2,
\end{multline}
where matrices $\M{A} \in \Real^{I \times R}, \M{B} \in \Real^{J \times R}$ and $\M{C} \in \Real^{K \times R}$
are the \emph{factor matrices} of $\T{X}$ extracted using a CANDECOMP/PARAFAC (CP) model \cite{CaCh70,Ha70,Hi27a}.
The CP model is one of the most commonly used tensor models in the literature (for a list of
CP applications, see reviews \cite{AcYe09, KoBa09}). Here, we use the notation $\T{X} = \KOp{\M{A},\M{B},\M{C}}$
to denote the CP model. Similarly, matrices $\M{A}$ and $\M{V} \in \Real^{M \times R}$ are the factor matrices
extracted from matrix $\M{Y}$ through matrix factorization. The formulation in \Eqn{fABCV} easily extends to multiple matrices
and tensors, e.g.,  as shown in \Fig{CoupledEx}. We also note that we focus on the least squares error in this paper, but our algorithms can be extended to other loss functions such as Bregman information metric used in, e.g., \cite{BaBaMe07}. We briefly illustrate two motivating applications of coupled matrix and tensor factorizations.

\begin{figure}[thb]
\centering
\includegraphics[width=2.5in]{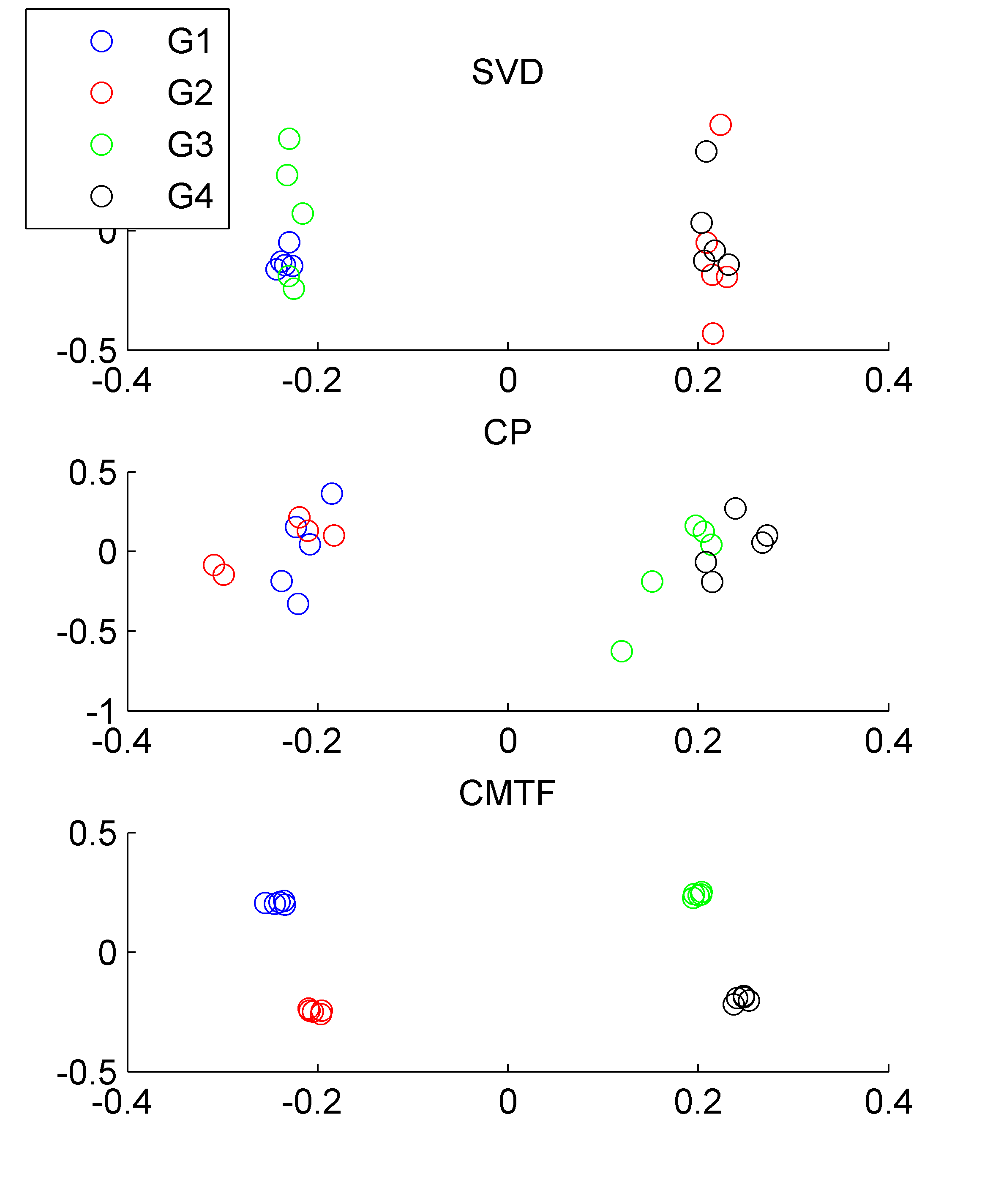}
\caption[Clustering Example]{Clustering based on matrix SVD factorization of $\M{Y}$ vs.\@
CP tensor factorization of $\T{X}$ vs.\@ coupled matrix-tensor factorization of $\T{X}$ and $\M{Y}$.
The subplots present the scatter plots showing the first factor plotted
against the second factor in the first mode.} \label{fig:Clustering}
\end{figure}

\noindent \textbf{Example 1: Clustering.} Joint analysis of data from multiple sources may capture fine-grained clusters that would not
be captured by the individual analysis of each data set. Suppose that there is a set of customers and there are two sources
of information about these customers, $\T{X} \in \Real^{I \times J \times K}$ and
$\M{Y} \in \Real^{I \times M}$, one storing information about which items
customers have bought over a period of time, and the other showing where customers live.
Within this set of customers, there are 4 groups: $\{G_1,G_2,G_3,G_4\}$
and each group consists of people who live in the same
neighborhood and have an interest in similar items.
Imagine that the matrix $\M{Y}$ can only discriminate between $(G_1 \cup G_3)$ and $(G_2 \cup G_4)$. A rank\nobreakdash-2 matrix SVD factorization of $\M{Y}$ would yield factors that could be used to cluster the customers as in the top plot shown in \Fig{Clustering} (SVD), failing to fully separate the four groups. Conversely, imagine that the tensor $\M{X}$ only has enough information to discriminate between $(G_1 \cup G_2)$ and $(G_3 \cup G_4)$. In this case, a rank-2 CP factorization of the tensor would still only separate the data into two groups, albeit two different groups, as illustrated in the middle plot in \Fig{Clustering} (CP). If, however, we jointly factor the matrix and tensor simultaneously using CMTF, the four groups are completely separated, as shown in the bottom plot of \Fig{Clustering} (CMTF). Details of the data generation for this example are provided in the appendix. \hfill$\square$

\begin{figure}[thb]
\centering
\includegraphics[width=2in]{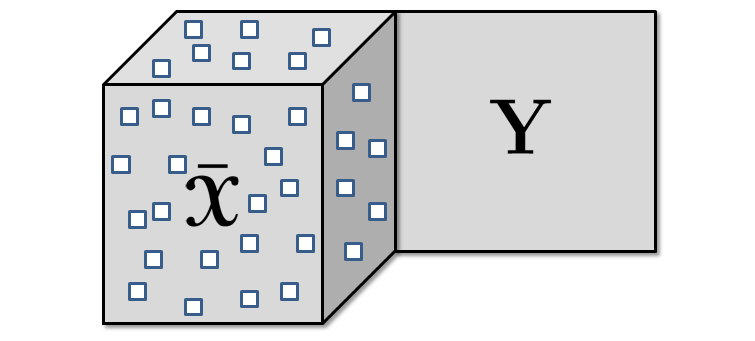}
\caption[Coupled Data Sets with Missing Entries]{Incomplete tensor $\Tbar{X}$ and
matrix $\M{Y}$ coupled in the first mode.} \label{fig:CoupledEx_Missing}
\end{figure}

\noindent \textbf{Example 2: Missing Data Recovery.} CMTF can be used for missing data recovery
when data from different sources have the same underlying low-rank structure (at least in one mode) but some
of the data sets have missing entries (\Fig{CoupledEx_Missing}). If a matrix or a higher-order tensor has a low-rank structure, it is possible to recover the missing entries using a limited number of data entries \cite{AcDuKoMo10,CaTa09}. However, if there is a large amount of missing data, then the analysis of a single data set is no longer enough for accurate data recovery. Here, we provide an example illustrating that missing entries
can still be recovered accurately using CMTF even when the analysis of a single data set fails to do so.

Suppose we have a tensor $\T{X}$ and a matrix $\M{Y}$ computed as $\T{X} = \KOp{\M{A},\M{B},\M{C}}$
and $\M{Y} = \M{A}\M{V}\Tra$, where matrices $\M{A} \in \Real ^{I \times R}$,
$\M{B} \in \Real^{J \times R}, \M{C} \in \Real^{K \times R}$ and $\M{V} \in \Real^{M \times R}$
are generated using random entries drawn from the standard normal distribution. $\Tbar{X} \in \Real^{I \times J \times K}$ is constructed
by randomly setting $M\%$ of the entries of $\T{X}$ to be missing, i.e., $\Tbar{X}=\T{W}*\T{X}$, where the binary tensor
$\T{W}$ is the same size as tensor $\T{X}$ and $\TE{W}{ijk}=0$ if we want to set $\TE{X}{ijk}$ to missing. In order to recover the missing entries,
one approach is to fit an $R$-component CP model to $\Tbar{X}$ and use the extracted
factors to recover missing entries. An alternative approach is to fit an $R$-component CMTF model to  $\Tbar{X}$
and $\M{Y}$ by extracting a common factor matrix in the first mode and then make use of CMTF factors to recover the missing entries.

\begin{figure}[t]
\centering
\includegraphics[width=2.5in]{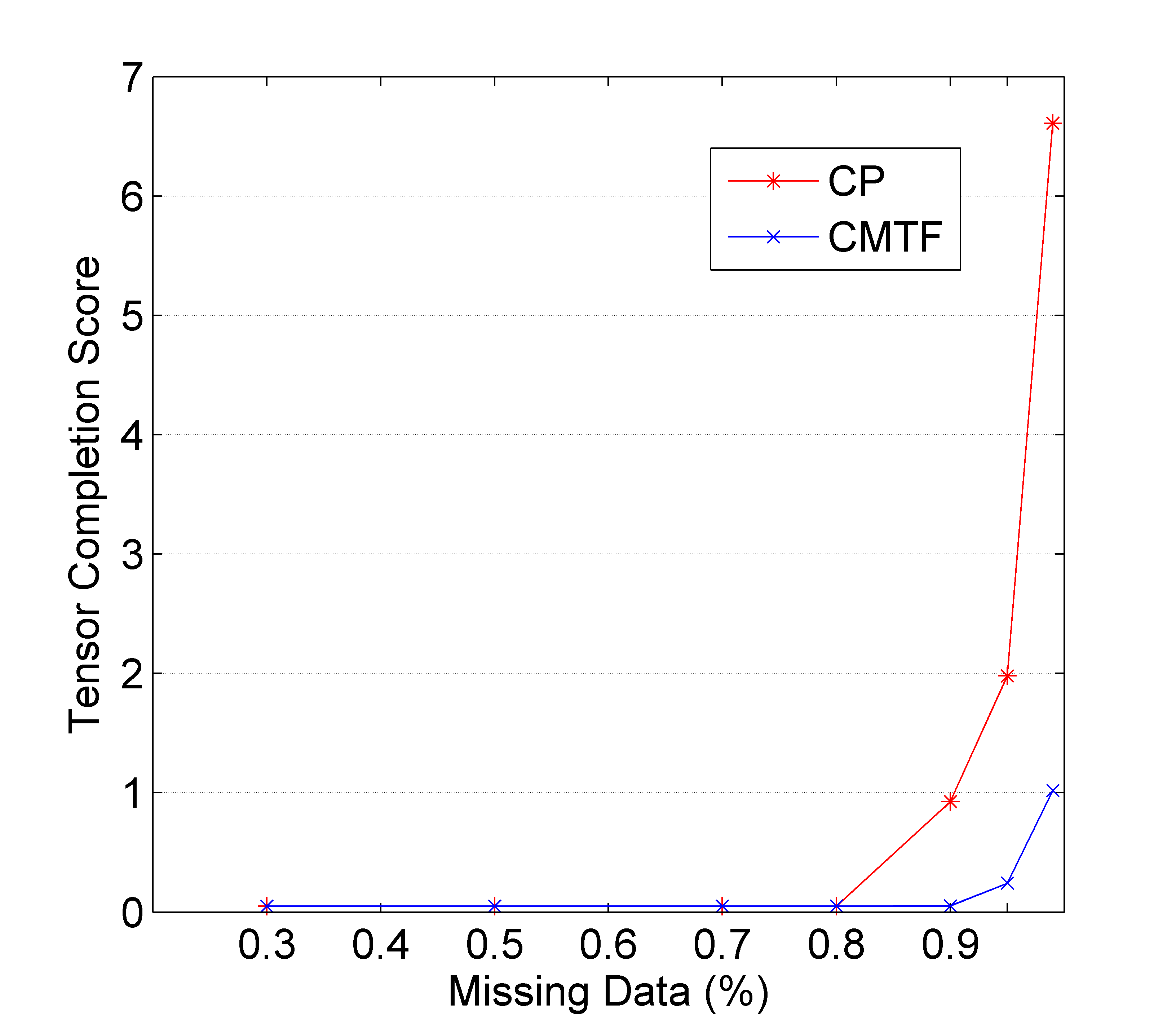}
\caption[Missing Data Recovery]{Missing data recovery using CP factorization of $\Tbar{X}$ vs.
CMTF of $\Tbar{X}$ and $\M{Y}$. Tensor Completion Score is the recovery error measure used here and is defined as $\frac{\norm{(1-\T{W})*(\T{X}-\That{X})}}{\norm{(1-\T{W})*\T{X}}}$,
where $\That{X}$ corresponds to the data reconstructed using the extracted factor matrices, i.e., $\Mhat{A},\Mhat{B},$ and $\Mhat{C}$, as $\That{X}=\KOp{\Mhat{A},\Mhat{B},\Mhat{C}}$.} \label{fig:Missing}
\end{figure}

\Fig{Missing} illustrates how the recovery error behaves for
different amounts of missing data. We observe that CP is accurate in terms of recovering missing entries
if less than $80\%$ of the entries are missing. However, there is a sharp in
error as we further increase the amount of missing entries. On the other hand, CMTF can compute factors with low recovery error
for higher amounts of missing data; only for problems with  more that $90\%$ missing data does the recovery error increase\footnote{Note that
the amount of missing data where the recovery error makes a sharp
increase may change depending on the values of $I,J,K,V$ and
$R$. For example, with small data sizes and large $R$, dealing with missing
data is challenging\cite{AcDuKoMo10}.}.\hfill$\square$

In order to solve coupled matrix and tensor factorization, various algorithms have been proposed in the literature using a variety of
loss functions \cite{BaBaMe07,LiSuCaKo09}. These algorithms all propose an alternating scheme,
where the basis for each entity type is determined
one at a time. In this paper, we focus solely on the $L_2$ loss function, which is often solved using alternating least squares (ALS).
Especially when fitting tensor models, ALS is the most commonly used algorithm due to its speed and ease of implementation.
On the other hand, ALS suffers from several problems: (i) It may fail to find the underlying components
accurately if the number of components is not correctly estimated \cite{AcDuKo11a,ToBr06}; (ii) In the presence
of missing data, ALS-based imputation techniques may suffer from poor convergence \cite{BuFi05} and
do not scale to large-scale data sets \cite{AcDuKoMo10}. Therefore, unlike the previous work using alternating
least squares algorithms, we propose an all-at-once optimization approach solving for all variables simultaneously.
Our contributions in this paper are
\begin{compactitem}
\item Developing an algorithm
called CMTF-OPT (CMTF-OPTimization) based on first-order
optimization to solve coupled matrix and tensor factorization problem for a given number of components ($R$).%
\footnote{Determining the number of components $R$ in CMTF remains as a challenge just like computing the tensor rank, which is NP-hard \cite{Ha90}.}
\item Extending
CMTF-OPT algorithm to handle incomplete data sets, i.e., data with
missing or unknown entries.
\item Demonstrating that CMTF-OPT is more accurate than an alternating
least squares approach using numerical experiments.
\end{compactitem}

This paper is organized as follows. In \Sec{notation}, we introduce the notation for tensors and tensor operations.
After discussing the related work on coupled data analysis in \Sec{related}, we introduce our CMTF-OPT algorithm and its extension
to data with missing entries in \Sec{algorithm}. \Section{experiments} describes numerical experiments and demonstrates
how CMTF-OPT compares with CMTF-ALS in terms of capturing the underlying factors in coupled data sets.
Finally, we conclude with future research directions in \Sec{conclusions}.


\section{Notation and Background}
\label{sec:notation}

Tensors of order $N \geq 3$ are denoted by Euler script letters
($\T{X},\T{Y},\T{Z}$), matrices are denoted by boldface capital letters
($\M{A},\M{B},\M{C}$), vectors are denoted by boldface lowercase letters
($\V{A},\V{B},\V{C}$), and scalars are denoted by lowercase letters ($a$, $b$, $c$).
Columns of a matrix are denoted by boldface lower letters with a
subscript, e.g., $\MC{A}{r}$ is the rth column of matrix $\M{A}$.
Entries of a matrix or a tensor are denoted by lowercase letters with
subscripts, i.e., the $(i_1, i_2, \dots, i_N)$ entry of an $N$-way
tensor $\T{X}$ is denoted by $\TE{X}{\SubscriptN{i}{N}}$.

Given a matrix $\M{A}$ of size $I \times J$, $\VC{\M{A}}$ stacks the columns of the matrix
and forms a vector of length $IJ$:
\begin{displaymath}
  \VC{\M{A}} =
  \begin{bmatrix}
    \\[-1em]
    \MC{A}{1} \\
    \vdots \\
    \MC{A}{J} \\[0.5em]
  \end{bmatrix}
  \in \Real^{IJ}.
\end{displaymath}

Given two matrices $\M{A} \in \Real^{I \times K}$ and $\M{B} \in \Real^{J \times K}$,
their Khatri-Rao product is denoted by $\M{A} \Khat \M{B}$ and defined as columnwise
Kronecker product. The result is a matrix of size
$(IJ) \times K$ and defined by
\begin{displaymath}
    \M{A} \Khat \M{B} =
    \begin{bmatrix}
        \MC{A}{1} \Kron \MC{B}{1}
        & \MC{A}{2} \Kron \MC{B}{2}
        & \cdots
        & \MC{A}{K} \Kron \MC{B}{K}
    \end{bmatrix},
\end{displaymath}
where $\Kron$ denotes Kronecker product. For more details on properties of Kronecker and Khatri-Rao
products, see \cite{KoBa09}.

An $N$-way tensor can be rearranged as a matrix; this is called \emph{matricization}.
The mode-$n$ matricization of a tensor $\T{X} \in \Real^{\SizeN{I}{N}}$ is denoted
by $\Mz{X}{n}$ and arranges the mode-$n$ one-dimensional ``fibers'' to be the columns of
the resulting matrix.

Given two tensors $\T{X}$ and $\T{Y}$ of equal size $\SizeN{I}{N}$,
their Hadamard (elementwise) product is denoted by $\T{X} \Hada \T{Y}$
and defined as
\begin{displaymath}
  \TEP{\T{X} \Hada \T{Y}}{\SubscriptN{i}{N}} =
  \TE{X}{\SubscriptN{i}{N}} \TE{Y}{\SubscriptN{i}{N}}
\end{displaymath}
for all $i_n \in \{1,\dots,I_n\}$ and $n \in \{1,\dots,N\}$. Their
inner product, denoted by $\ip{\T{X}}{\T{Y}}$, is the sum of the products 
of their entries, i.e.,
\begin{displaymath}
\ip{\T{X}}{\T{Y}}
=
\SumN{i}{I}{N}
\TE{X}{\SubscriptN{i}{N}}
\TE{Y}{\SubscriptN{i}{N}}.
\end{displaymath}
For a tensor $\T{X}$ of size $\SizeN{I}{N}$, its \emph{norm} is
\begin{inlinemath}
  \norm{ \T{X} } = \sqrt{\ip{\T{X}}{\T{X}}}.
\end{inlinemath}
For matrices and vectors, $\|\cdot\|$ refers to the analogous
Frobenius and two-norm, respectively.

Given a sequence of matrices $\Mn{A}{n}$ of size $I_n \times R$ for
$n=1,\dots,N$, the notation
$\KOp{\Mn{A}{1},\Mn{A}{2},\dots,\Mn{A}{N}}$ defines an $\SizeN{I}{N}$
tensor whose elements are given by
\begin{displaymath}
  \TEP{\KOp{\Mn{A}{1},\Mn{A}{2},\dots,\Mn{A}{N}}}{\SubscriptN{i}{N}}
  = \sum_{r=1}^R \prod_{n=1}^N \MnE{a}{n}{i_nr},
\end{displaymath}
for $i_n \in \{1,\dots,I_n\}, n \in \{1,\dots,N\}$. For just two matrices,
this reduces to
\begin{inlinemath}
  \KOp{\M{A},\M{B}} = \M{A}\M{B}\Tra.
\end{inlinemath}

\section{Related Work in Data Fusion}
\label{sec:related}

Data fusion, also called collective data analysis, multi-block, multi-view
or multi-set data analysis, has been a topic of interest
in different fields for decades. First, we briefly discuss data fusion
techniques for multiple data sets each represented as a matrix and then focus on techniques proposed for coupled
analysis of heterogenous data sets.

\subsection{Collective Factorization of Matrices}
The analysis of data from multiple sources attracted considerable attention
in the data mining community during the Netflix Prize competition \cite{KoBeVo09},
where the goal was to make accurate predictions about movie ratings.
In order to achieve better rating predictions, additional data sources complementing
user ratings such as tagging information have been exploited; e.g., users tag movies \cite{WaWaWu10}
as well as features of movies such as movie types or movie players .
Singh and Gordon \cite{SiGo08a} proposed Collective Matrix Factorization (CMF) to
take advantage of correlations between different data sets and simultaneously factorize
coupled matrices. Gi-ven two matrices $\M{X}$ and $\M{Y}$ of
size $I \times M$ and $I \times L$, respectively, CMF can be formulated as
\begin{equation}
  \label{eq:M}
    f(\M{U},\M{V},\M{W}) = \norm{\M{X} - \M{U}\M{V} \Tra}^2 + \norm{\M {Y}-\M{U}\M{W}\Tra}^2,
\end{equation}
where $\M{U},\M{V}$ and $\M{W}$ are factor matrices of size $I \times R$,
$M \times R$ and $L \times R$, respectively and $R$ is the number of factors.
This formulation is a special case of the general approach introduced
in \cite{SiGo08a}, which extends to different loss functions. Earlier, Long et al.
\cite{LoZhWuYu06,LoWuZhYu06} had also studied collective matrix factorization
using a different matrix factorization scheme than $\M{X}=\M{U}\M{V}\Tra$.
The proposed approaches in those studies are based on alternating algorithms solving
the collective factorization problem for one factor matrix at a time.

Analysis of multiple matrices dates back to one of the earliest models aiming
to capture the common variation in two data sets, i.e., Canonical Correlation
Analysis (CCA) \cite{Ho36}. Later, other studies followed CCA by extending it to more
than two data sets \cite{Ke69}, focusing on simultaneous factorization of Gramian matrices
\cite{Le66} and working on PCA of multiple matrices \cite{KaHuSiPaWa00,WeKoMa98}.
Moreover, several approaches for simultaneous factor analysis have been developed for 
specific applications as well; e.g., population differentiation in biology \cite{Th88}, blind source 
separation \cite{ZiLaNoMu04}, multimicrophone speech filtering \cite{DoMo02}, and microarray 
data analysis \cite{AlBrBo03,Ba07,Ba08}.

Tensor factorizations \cite{AcYe09,KoBa09,SmBrGe04} can also be considered as one way of
analyzing multiple matrices. For instance, when a tensor model is fit to a third-order tensor,
multiple coupled matrices are analyzed simultaneously. Nevertheless, neither CMF nor
tensor factorizations can handle coupled analysis of heterogeneous data, which we address next.

\subsection{Collective Factorization of Mixed Data}
As described in \Sec{introduction}, heterogeneous data consists of data sets of
different orders, i.e., both matrices and higher-order tensors. The formulation in \Eqn{M}
can be extended to heterogeneous data sets: Given a tensor $\T{X}$ and a matrix $\M{Y}$ of sizes
$I \times J \times K$ and $I \times M$, respectively, we can formulate their factorization coupled
in the first mode as shown in \Eqn{fABCV}.
Note that \Eqn{fABCV} can be considered as a special case of the approach
introduced by Smilde et~al.~\cite{SmWeBo00} for multi-way multi-block data analysis, where
different tensor models can be fit to higher-order data sets and the factor matrices corresponding to the
coupled modes do not necessarily match. The same formulation as in \Eqn{fABCV}
has recently been studied in psychometrics as Linked-Mode PARAFAC-PCA \cite{Wilderjans2009}. As a more general framework,
not restricted to squared Euclidean distance, Banerjee et~al.~\cite{BaBaMe07} introduced a multi-way clustering approach
for relational and multi-relational data where coupled analysis of multiple data sets including higher-order
data sets were studied using minimum Bregman information. The paper \cite{LiSuCaKo09} also discussed coupled analysis
of multiple tensors and matrices using nonnegative factorization by formulating the problem using KL-divergence.
All these studies propose algorithms that are based on alternating approaches.

In this paper we focus on the squared Euclidean distance as the loss function as in \Eqn{fABCV}. Rewriting $\T{X} =
\KOp{\M{A},\M{B},\M{C}}$ as $\TM{X}{1}=\M{A} (\M{C}\Khat \M{B})\Tra$,  $\TM{X}{2}=\M{B}
(\M{C}\Khat \M{A})\Tra$, etc., for the different modes of $\T{X}$, we summarize the steps of an alternating least squares
algorithm for solving \Eqn{fABCV} in \Fig{ALS}. The main loop in the algorithm is often terminated as a function of the relative
change in function value (e.g., as defined below in \Eqn{relfunc}), a function of the relative change in factor matrices, and/or after a prescribed number of iterations.

\begin{figure}[t!]
\centering
\begin{boxedminipage}{3.3in}
\begin{algorithmic}
  \WHILE {not ``converged''}
    \STATE rescale all factors to unit Frobenious norm
    \STATE {solve for $\M{A}$ (for fixed $\M{B},\M{C},\M{V}$)\\
     $\quad \min_{\M{A}}  \norm{ \left[\TM{X}{1} \quad \M{Y}\right] - \M{A} \left[(\M{C}\Khat \M{B})\Tra \quad \M{V}\Tra\right]}^2$ }
    \STATE {solve for $\M{B}$ (for fixed $\M{A},\M{C},\M{V}$) \\
     $\quad \min_{\M{B}}  \norm{ \TM{X}{2}  - \M{B} (\M{C}\Khat \M{A})\Tra}^2$}
    \STATE {solve for $\M{C}$ (for fixed $\M{A},\M{B},\M{V}$)\\
     $\quad \min_{\M{C}}  \norm{ \TM{X}{3}  - \M{C} (\M{B}\Khat \M{A})\Tra}^2$}
    \STATE {solve for $\M{V}$ (for fixed $\M{A},\M{B},\M{C}$)\\
     $\quad \min_{\M{V}}  \norm{ \M{Y}  - \M{A}\M{V}\Tra}^2$}
 \ENDWHILE
\end{algorithmic}
\end{boxedminipage}
\caption[CMTF-ALS]{CMTF-ALS: Alternating Least Squares Algorithm for coupled matrix
and tensor factorization of a third-order tensor $\T{X}$ and matrix $\M{Y}$ coupled in the first mode.} \label{fig:ALS}
\end{figure}

ALS-based algorithms are simple to implement and computationally efficient;
however, ALS has shown to be error-prone when fitting a CP model in the case of
overfactoring, i.e., the number of extracted components
is more than the true number of underlying components \cite{AcDuKo11a, ToBr06}. Furthermore, in the presence of
missing data, ALS-based techniques may have poor convergence \cite{BuFi05}
and do not scale to very large data sets \cite{AcDuKoMo10}. On the other hand,
all-at-once optimization, in other words, solving for all CP factor matrices simultaneously
has shown to be more robust to overfactoring \cite{AcDuKo11a, ToBr06} and easily extends to handle
data with missing entries even for very large data sets \cite{AcDuKoMo10}. Therefore, in order to
deal with these issues, we develop an algorithm called CMTF-OPT, which formulates coupled analysis of heterogeneous
data sets just like \cite{SmWeBo00,Wilderjans2009} but solves for all
factor matrices simultaneously using a gradient-based optimization approach.


\section{CMTF-OPT Algorithm}
\label{sec:algorithm}
In this section, we consider joint analysis of a matrix and an $N$th-order tensor with one mode in common, where the tensor is factorized using an $R$-component CP model and the matrix is factorized by extracting $R$ factors using matrix factorization. Let $\T{X} \in \Real^{\SizeN{I}{N}}$ and $\M{Y} \in \Real^{I_1 \times M}$ have the $n{th}$ mode in common, where $n \in \{1,\dots,N\}$. The objective function for coupled analysis of these two data sets is defined by
\begin{multline}
\label{eq:fNway}
   f(\Mn{A}{1},\Mn{A}{2},\dots,\Mn{A}{N},\M{V})\\
     =  \frac{1}{2}\norm{\T{X} - \KOp{\Mn{A}{1},\dots,\Mn{A}{N}}}^2 + \frac{1}{2}\norm{\M{Y} - \Mn{A}{n}\M{V}\Tra}^2
\end{multline}
Our goal is to find the matrices $\Mn{A}{i} \in \Real^{I_i \times R}$ for $i=1, 2, ...N$ and matrix $\M{V} \in \Real^{M \times R}$ that minimize the objective in \Eqn{fNway}. In order to solve this optimization problem, we can compute the gradient and then use any first-order optimization algorithm \cite{NoWr06}. Next, we discuss the computation of the gradient for \Eqn{fNway}.

We can rewrite \Eqn{fNway} as two components, $f_1$ and $f_2$:
\begin{multline*}
   f = \frac{1}{2} \underbrace{\norm{\T{X} - \KOp{\Mn{A}{1},\dots,\Mn{A}{N}}}^2}_{f_1(\Mn{A}{1},\Mn{A}{2},\dots,\Mn{A}{N})} + \frac{1}{2} \underbrace{\norm{\M{Y} - \Mn{A}{n}\M{V}\Tra}^2}_{f_2(\Mn{A}{n},\M{V})}
\end{multline*}

The partial derivative of ${f_1}$ with respect to ${\Mn{A}{i}}$ has been derived in \cite{AcDuKo11a} so we just present the results here.

Let $\T{Z}=\KOp{\Mn{A}{1},\dots,\Mn{A}{N}}$, then
\begin{displaymath}
  \FD{f_1}{\Mn{A}{i}} = (\Mz{Z}{i}-\Mz{X}{i})\Mn{A}{-i}
\end{displaymath}
where
\begin{displaymath}
\Mn{A}{-i} = \Mn{A}{N} \Khat \cdots \Khat \Mn{A}{i+1} \Khat
\Mn{A}{i-1} \Khat \cdots \Khat \Mn{A}{1},
\end{displaymath}
for $i=1,\dots,N.$

The partial derivatives of the second component, ${f_2}$, with respect to ${\Mn{A}{i}}$  and $\M{V}$ can be computed as
\begin{align*}
  \FD{f_2}{\Mn{A}{i}} & =
   \begin{cases}
     -\M{Y}\M{V} + \Mn{A}{-i}\M{V}\Tra\M{V}, & \text{for $i=n$},\\
     0 & \text{for $i \neq n $},\\
  \end{cases}\\
  \FD{f_2}{\M{V}} & = -\M{Y}\Tra \Mn{A}{i} + \M{V}{\Mn{A}{i}}\Tra\Mn{A}{i} .
\end{align*}

Combining the above results, we can compute the partial derivative of $f$ with respect to factor matrix $\Mn{A}{i}$, for $i=1,2, ...,N$, and $\M{V}$ as:
\begin{align*}
  \FD{f}{\Mn{A}{i}} & =
    \FD{f_1}{\Mn{A}{i}} + \FD{f_2}{\Mn{A}{i}} \\
  \FD{f}{\M{V}} & = \FD{f_2}{\M{V}}
\end{align*}

Finally, the gradient of $f$, which is a vector of size $P=R(\sum_{n=1}^N I_n + M)$,
can be formed by vectorizing the partial derivatives with respect to
each factor matrix and concatenating them all, i.e.,
\begin{displaymath}
 \nabla f=
   \begin{bmatrix}
    \\[-1em]
    \VC{\FD{f}{\Mn{A}{1}}} \\
    \vdots \\
    \VC{\FD{f}{\Mn{A}{N}}} \\[0.5em]
    \VC{\FD{f}{\M{V}}} \\[0.5em]
  \end{bmatrix}
\end{displaymath}

\subsection{CMTF-OPT for Incomplete Data}
\label{sec:algorithm_missing}
In the presence of missing data, it is still possible
to do coupled analysis by ignoring the missing
entries and fitting the tensor and/or the matrix model to the known data entries.
Here we study the case where tensor $\T{X}$ has missing entries as in \Fig{CoupledEx_Missing}. Let
$\T{W} \in \Real^{\SizeN{I}{N}}$ indicate the missing entries of $\T{X}$ such that
\begin{displaymath}
  \label{eq:W}
  \TE{W}{\SubscriptN{i}{N}} =
  \begin{cases}
    1 & \text{if $\TE{X}{\SubscriptN{i}{N}}$ is known},\\
    0 & \text{if $\TE{X}{\SubscriptN{i}{N}}$ is missing},
  \end{cases}
\end{displaymath}
for all  $i_n \in \{1,\dots,I_n\}$ and $n \in \{1,\dots,N\}$. We can then modify the
objective function \Eqn{fNway} as
\begin{multline*}
   f_{\T{W}}(\Mn{A}{1},\Mn{A}{2},\dots,\Mn{A}{N},\M{V})\\
     =  \frac{1}{2}\underbrace{\norm{\T{W} \Hada \left(\T{X} - \KOp{\Mn{A}{1},\dots,\Mn{A}{N}}\right)}^2}_{f_{\T{W}_1}(\Mn{A}{1},\Mn{A}{2},\dots,\Mn{A}{N})} + \frac{1}{2}\norm{\M{Y} - \Mn{A}{n}\M{V}\Tra}^2
\end{multline*}

The first term, $f_{\T{W}_1}$, corresponds to the weighted least squares problem for fitting a CP model while the second term stays the same
as in \Eqn{fNway}. The partial derivative of $f_{\T{W}_1}$ with respect to $\Mn{A}{i}$ can be computed as in \cite{AcDuKoMo10} as
\begin{displaymath}
  \FD{f_{\T{W}_1}}{\Mn{A}{i}} = \left(\Mz{W}{i} \Hada \Mz{Z}{i}-\Mz{W}{i} \Hada \Mz{X}{i}\right)\Mn{A}{-i},
\end{displaymath}
for $i=1,\dots,N.$ Since the partial derivatives of the second term do not change, we can write the partials of $f$ with respect to $\Mn{A}{i}$
for $i=1,2, ...,N$, and $\M{V}$ as:
\begin{align*}
  \FD{f_{\T{W}}}{\Mn{A}{i}} & =
  \begin{cases}
    \FD{f_{\T{W}_1}}{\Mn{A}{i}} & \text{for $i \in \{1,\dots,N\}/\{n\}$}, \\
    \FD{f_{\T{W}_1}}{\Mn{A}{i}} + \FD{f_2}{\Mn{A}{i}} & \text{for $i =n $}. \\
  \end{cases}\\
  \FD{f_{\T{W}}}{\M{V}} & = \FD{f_2}{\M{V}}
\end{align*}
The gradient of $f_{\T{W}}$, $\nabla f_{\T{W}}$, is also a vector of size $P=R(\sum_{n=1}^N I_n + M)$ and can be formed in the same way as $\nabla f$.

Once we have the function, $f$ (or $f_{\T{W}}$), and gradient, $\nabla f$ (or $\nabla f_{\T{W}}$), values, we can use any gradient-based optimization algorithm to compute the factor matrices. For the results presented in this paper, we use the Nonlinear Conjugate Gradient (NCG) with Hestenes-Steifel updates \cite{NoWr06} and the
Mor\'{e}-Thuente line search \cite{MoTh94} as implemented in the Poblano Toolbox \cite{DuKoAc10}.


\section{Experiments}
\label{sec:experiments}
We compare the performance of the proposed CMTF-OPT algorithm with the
ALS-based approach (i.e., CMTF-ALS as shown in \Fig{ALS}) in terms of accuracy
using randomly generated matrices and tensors. Our goal is to see whether
the algorithms can capture the underlying factors in the data (i) when  $\bar R = R$
factors are extracted, and (ii) in the case of overfactoring, e.g., when $\bar R = R+1$
factors are extracted, where $\bar R$ and $R$ denote the extracted and true number of factors.

\subsection{Data Generation}
In our experiments, three different scenarios are used (see \Tab{ACC1}). In the first case, we have a tensor $\T{X}$
and a matrix $\M{Y}$ coupled in the first mode. In the second case, we have two third-order tensors $\T{X}$ and $\T{Y}$ coupled in one dimension. The third case has a third-order tensor $\T{X}$ and two matrices $\M{Y}$ and $\M{Z}$ such that the tensor shares one mode with each one of the matrices.

To construct the tensors and matrices in each scenario, random factor matrices with entries following standard normal distribution
are generated and a tensor (or a matrix) is formed based on a CP (or a matrix) model. For example, for the first scenario, we generate random
factor matrices $\M{A} \in \Real^{I \times R}$, $\M{B} \in \Real^{J \times R}$ and $\M{C} \in \Real^{K \times R}$,
and form a third-order tensor $\T{X} \in \Real^{I \times J \times K}$ based on a CP model. Gaussian noise is later added to the tensor, i.e., $\T{X}=\KOp{\M{A}, \M{B},\M{C}} + \eta \T{N} \frac{\norm{\KOp{\M{A}, \M{B},\M{C}}}}{\norm{\T{N}}}$, where $\T{N} \in \Real^{I \times J \times K}$ corresponds to the random noise tensor and $\eta$ is used to adjust the noise level. Similarly, we generate a factor matrix $\M{V} \in \Real^{I \times R}$ and form matrix $\M{Y}=\M{A}\M{V}\Tra+\eta\M{N} \frac{\norm{\M{A}\M{V}\Tra}}{\norm{\M{N}}}$. The matrix $\M{N} \in \Real^{I \times R}$ is a matrix with entries drawn from the standard normal distribution and is used to introduce differing amounts of noise in the data. In our experiments, we use three different noise levels, i.e., $\eta=0.1, 0.25$ and $0.35$, and the true number of factors is $R=3$.

\subsection{Performance Metric}
We fit an $\bar R$-component CMTF model using CMTF-OPT and CMTF-ALS algorithms, where $\bar R=R$ and $\bar R=R+1$, and compare the algorithms in terms of accuracy; in other words, how well the factors extracted by the algorithms match with the original factors used  to generate the data. For instance, for the first scenario described above, let $\M{A}, \M{B}, \M{C}$ and $\M{V}$ be the original factor matrices and  $\Mhat{A}, \Mhat{B}, \Mhat{C}$ and $\Mhat{V}$ be the extracted factor matrices. We quantify how well the extracted factors match with the original ones using the factor match score (FMS) defined as follows:
\begin{multline}
\label{eq:FMS}
  \text{FMS} = \min_{r}(1-\frac{|\xi_r - \hat\xi_r |}{\max(\xi_r,\hat\xi_r)}) |\MCTra{A}{r}\MhatC{A}{r}\MCTra{B}{r}\MhatC{B}{r}\MCTra{C}{r}\MhatC{C}{r}\MCTra{V}{r}\MhatC{V}{r}|.
\end{multline}
Here, the columns of factor matrices are normalized to unit norm and $\xi_r$ denotes the weight for each factor $r$, for $r=1,2,...R$. The weight for each factor is computed as follows: We can rewrite $\T{X}=\KOp{\M{A}, \M{B},\M{C}}$ as $\T{X}=\sum_{r=1}^R\lambda_r \MC{A}{r} \circ \MC{B}{r} \circ \MC{C}{r}$, where $\circ$ denotes $n$-mode vector outer product; the columns of the factor matrices are normalized to unit norm and $\lambda_r$ is the product of the vector norms in each mode. Similarly, we can rewrite $\M{Y}=\M{A}\M{V}\Tra$ as $\M{Y}=\sum_{r=1}^R\alpha_r \MC{A}{r} \circ \MC{V}{r}$ by normalizing the matrix columns and using $\alpha_r$ to denote the product of vector norms. We define the norm of each CMTF component, i.e., $\xi_r$, as $\xi_r=\lambda_r + \alpha_r$. In \Eqn{FMS}, $\xi_r$ and $\hat\xi_r$ indicate the weights corresponding to the original and extracted $r{th}$ CMTF component, respectively. If the extracted and original factors perfectly match, FMS value will be 1. It is also considered a success if FMS is above a certain threshold, i.e., $(0.99)^N$, where $N$ is the number of factor matrices.

\subsection{Stopping Conditions}
As a stopping condition, both algorithms use the relative change in function value and stop when
\begin{equation}
\label{eq:relfunc}
\frac{|f_{\text{current}}-f_{\text{previous}}|}{f_{\text{previous}}} \leq 10^{-8},
\end{equation}
where $f$ is as given in \Eqn{fNway}. Additionally, for CMTF-ALS, the maximum number of iterations
is set to $10^4$. For CMTF-OPT, the maximum number of function values (which corresponds the number
of iterations in an ALS algorithm) is set to $10^4$ and the maximum number of iterations is set to
$10^3$. CMTF-OPT also uses the two-norm of the gradient divided by the number of entries in the gradient
and the tolerance for that is set to $10^{-8}$. In the experiments, the algorithms have generally stopped
due to the relative change in function value criterion. However, we also rarely observe that both CMTF-OPT and CMTF-ALS
have stopped since the algorithms have reached the maximum number of iterations, i.e., approximately $2\%$ of all runs for CMTF-OPT and $0.3\%$ of all runs for CMTF-ALS. Furthermore, for around $2\%$ of all runs, CMTF-OPT has stopped by satisfying the condition on the gradient.

\begin{table*}[ht!]
\footnotesize
\centering
\caption{Comparison of CMTF-OPT and CMTF-ALS for the cases where $\lambda_r=\alpha_r=...=1$, for $r=1,2,...R$.}
\small
\begin{tabular}{ccc|c c c |c c c|c c c|}
\cline{4-12}
    & & & & & & &  & & & & \\
    & & & \multicolumn{3}{c|} {\includegraphics[width =1.1in]{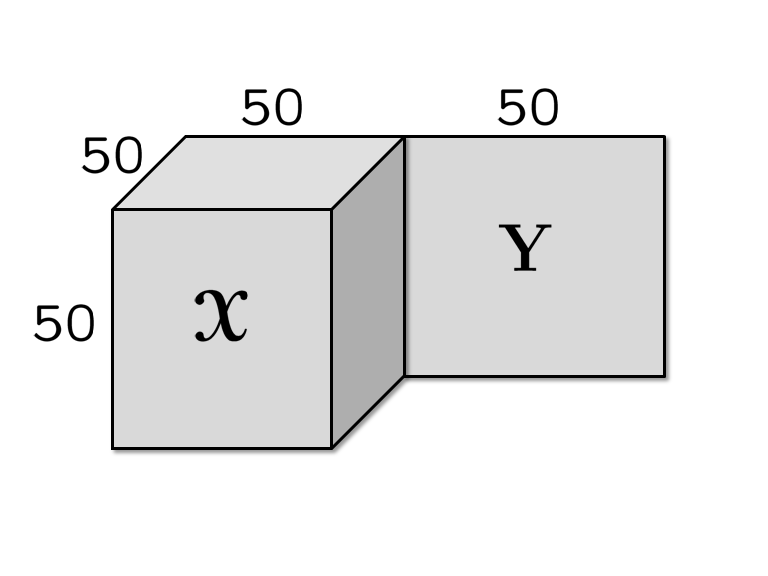}} & \multicolumn{3}{c|} {\includegraphics[width=1.1in]{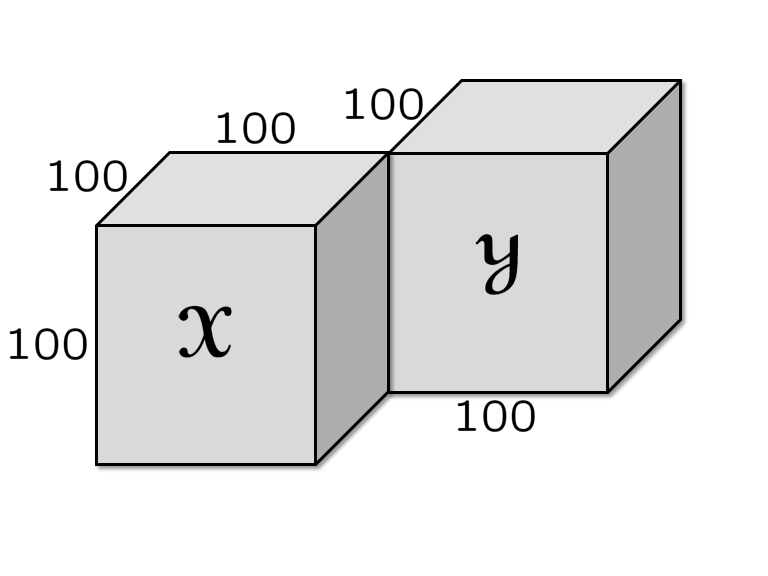}} & \multicolumn{3}{c|} {\includegraphics[width=1.1in]{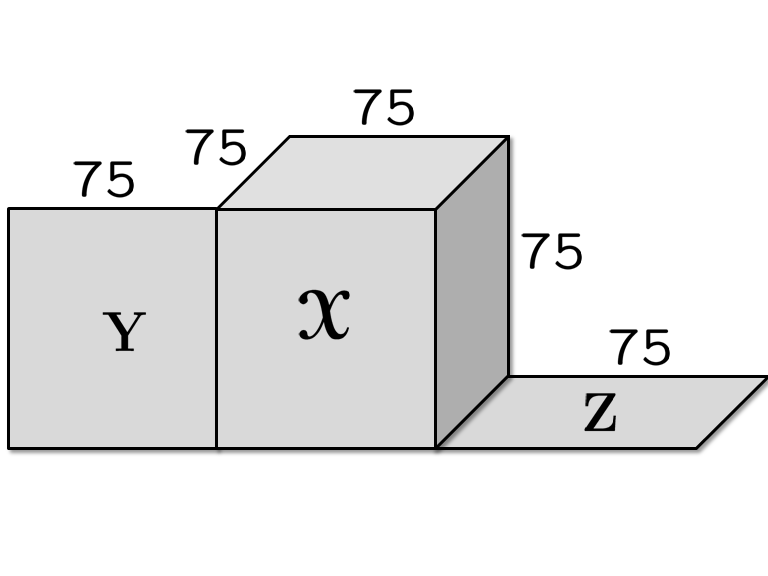}}\\
    & & & & & & &  & & & & \\ \hline
    \multicolumn{1}{|c|} {\multirow{2}{*}{\bf Noise}} & \multicolumn{1}{c|}{\multirow{2}{*}{\bf $\bar R$}} & \multirow{2}{*}{\bf ALG.} & \multicolumn{1}{c}{\bf Success} & \multicolumn{1}{c}{\bf Mean} & \multicolumn{1}{c}{\bf p-val} & \multicolumn{1}{|c}{\bf Success} & \multicolumn{1}{c}{\bf Mean } & \multicolumn{1}{c|}{\bf p-val} & \multicolumn{1}{c}{\bf Success} & \multicolumn{1}{c}{\bf Mean} & \multicolumn{1}{c|}{\bf p-val} \\
    \multicolumn{1}{|c|}{} & \multicolumn{1}{c|}{} & & \multicolumn{1}{c}{\bf (\%)} & \multicolumn{1}{c}{\bf FMS} & \multicolumn{1}{c}{} & \multicolumn{1}{|c}{\bf (\%)} & \multicolumn{1}{c}{\bf FMS} & \multicolumn{1}{c|}{} & \multicolumn{1}{c}{\bf (\%)} & \multicolumn{1}{c}{\bf FMS} & \multicolumn{1}{c|}{} \\ \hline \hline
    \multicolumn{1}{|c|} {\multirow{4}{*}{\bf$\eta$ = 0.10}} &  \multicolumn{1}{c|}{\multirow{2}{*}{\bf $R$}}  & {\bf OPT} & 100.0 & 1.00 & 9.3e-1 & 96.7 & 0.97 & 1.0 & 100.0 & 1.00 & 3.3e-1\\
    \multicolumn{1}{|c|} {} & \multicolumn{1}{c|}{} & {\bf ALS} & 100.0 & 1.00 & & 96.7 & 0.97 & & 96.7 & 0.96 &  \\ \cline{2-12}
    \multicolumn{1}{|c|} {} & \multicolumn{1}{c|}{\multirow{2}{*}{\bf $ R+1$}}  & {\bf OPT} & 96.7 & 0.97 & 3.3e-13 & 100.0 & 1.00 & 1.6e-8 & 96.7 & 0.97 & 3.8e-24 \\
    \multicolumn{1}{|c|} {} & \multicolumn{1}{c|}{} & {\bf ALS} & 3.3 & 0.29 &  & 6.7 & 0.71 & & 0.0 & 0.06 & \\  \hline \hline
    \multicolumn{1}{|c|} {\multirow{4}{*}{\bf$\eta$ = 0.25}} & \multicolumn{1}{c|}{\multirow{2}{*}{\bf $R$}} & {\bf OPT} & 100.0 & 0.99 & 4.5e-1 & 100.0 & 1.00 & 0.2e-1 & 96.7 & 0.96  & 3.3e-1 \\
    \multicolumn{1}{|c|} {} & \multicolumn{1}{c|}{} & {\bf ALS} & 100.0 & 0.99 &  & 100.0 & 1.00  & & 100.0 & 0.99 & \\ \cline{2-12}
    \multicolumn{1}{|c|} {} & \multicolumn{1}{c|}{\multirow{2}{*}{\bf $R+1$}}  & {\bf OPT} & 100.0 & 0.99 & 2.6e-11 & 100.0  & 1.00 & 1.8e-7 & 100.0 & 0.99 & 7.0e-16\\
    \multicolumn{1}{|c|} {} & \multicolumn{1}{c|}{} & {\bf ALS} & 6.7 & 0.34 &  & 16.7  & 0.70 & & 10.0 & 0.14 & \\  \hline \hline
    \multicolumn{1}{|c|} {\multirow{4}{*}{\bf $\eta$ = 0.35}} & \multicolumn{1}{c|}{\multirow{2}{*}{\bf $R$}} & {\bf OPT} & 100.0 & 0.99 & 5.8e-1 & 100.0 & 1.00 & 3.3e-1 & 100.0  & 0.99  & 3.6e-1\\
    \multicolumn{1}{|c|} {} & \multicolumn{1}{c|}{} & {\bf ALS} & 100.0 & 0.99 &  & 96.7  & 0.97  & & 100.0 & 0.99 & \\ \cline{2-12}
    \multicolumn{1}{|c|} {} & \multicolumn{1}{c|}{\multirow{2}{*}{\bf $R+1$}}  & {\bf OPT} & 90.0 & 0.92 & 1.8e-5 & 100.0 & 1.00 & 1.0e-7 & 86.7 & 0.88 & 2.9e-5\\
    \multicolumn{1}{|c|} {} & \multicolumn{1}{c|}{} & {\bf ALS} & 33.3 & 0.55 &  & 13.3 & 0.71 & & 36.7 & 0.39 & \\  \hline
\end{tabular}
\label{tab:ACC1}
\end{table*}

\subsection{Results} We demonstrate the performance of the algorithms in terms of accuracy in \Tab{ACC1} and \Tab{ACC2}.
\Tab{ACC1} presents the results for the experiments where factor matrices are generated with columns normalized to unit norm. For instance, for the example discussed above, this corresponds to generating tensor $\T{X}$ and matrix $\M{Y}$ using $\lambda_r=\alpha_r=1$, for $r=1,2,..R$. In \Tab{ACC1}, for all different scenarios, we observe that when the correct number of underlying factors are extracted, i.e., $\bar R = R$, the success ratios of both algorithms are compatible. Since we repeat our experiments with 30 different sets of factor matrices for each set of parameters, we report the average factor match score for 30 runs. The average scores for both algorithms are quite close and p\nobreakdash-values computed for paired-sample $t$-tests indicate that differences in the scores are not statistically significant. On the other hand, for all scenarios, when we look at the cases where the data is overfactored, i.e., $\bar R = R+1$, CMTF-OPT significantly outperforms CMTF-ALS. While CMTF-OPT algorithm is quite accurate in terms of recovering the underlying factors in the case of overfactoring, the accuracy of CMTF-ALS is quite low. We also observe that the increase in the noise level barely affects the accuracies of the algorithms.

In \Tab{ACC2}, we present the results for a harder set of experiments, where the norm of each tensor and matrix component, i.e., $\lambda_r,\alpha_r$ for $r=1,2,..R$, is a randomly assigned integer greater than or equal to 1.\footnote{Each norm is chosen as the absolute value of a number randomly chosen from $\mathcal{N}(0,25)$ rounded to the nearest integer plus 1.} Similar to the previous case, we observe that CMTF-OPT is more robust to overfactoring compared to CMTF-ALS. However, the accuracies reported in this table are lower compared to \Tab{ACC1}. In particular, as the noise level increases, it becomes harder to find the underlying components and accuracies drop even if the true number of underlying factors are extracted from the data. For instance, for the first scenario, when the noise level is 0.35, the accuracy of CMTF-OPT is $60\%$ for the case we extract the correct number of components. This is in part due to not being able to find $\xi_r$, for $r=1,2,..R$, accurately. If we slightly change the FMS such that we only require $\min_{r}|\MCTra{A}{r}\MhatC{A}{r}\MCTra{B}{r}\MhatC{B}{r}\MCTra{C}{r}\MhatC{C}{r}\MCTra{V}{r}\MhatC{V}{r}|$ to be greater than the threshold, then the average accuracy goes up to around $73\%$. The rest of the failing runs is due to the fact that extracted factors are distorted compared to the original factors used to generate the data. Note that factor match scores are still high which indicates that the factors are only slightly distorted.

\begin{table*}[ht!]
\centering
\caption{Comparison of CMTF-OPT and CMTF-ALS for the cases where  $\lambda_r,\alpha_r,... \geq 1$ for $r=1,2,...R$.}
\small
\begin{tabular}{ccc|c c c |c c c |c c c|}
\cline{4-12}
    & & & & & & &  & & & & \\
    & & & \multicolumn{3}{c|} {\includegraphics[width =1.1in]{exp1}} & \multicolumn{3}{c|} {\includegraphics[width=1.1in]{exp3}} & \multicolumn{3}{c|} {\includegraphics[width=1.1in]{exp4}}\\
    & & & & & & &  & & & & \\ \hline
    \multicolumn{1}{|c|} {\multirow{2}{*}{\bf Noise}} & \multicolumn{1}{c|}{\multirow{2}{*}{\bf $\bar R$}} & \multirow{2}{*}{\bf ALG.} & \multicolumn{1}{c}{\bf Success)} & \multicolumn{1}{c}{\bf Mean} & \multicolumn{1}{c}{\bf p-val} & \multicolumn{1}{|c}{\bf Success } & \multicolumn{1}{c}{\bf Mean } & \multicolumn{1}{c|}{\bf p-val} & \multicolumn{1}{c}{\bf Success } & \multicolumn{1}{c}{\bf Mean } & \multicolumn{1}{c|}{\bf p-val} \\
    \multicolumn{1}{|c|}{} & \multicolumn{1}{c|}{} & & \multicolumn{1}{c}{\bf (\%)} & \multicolumn{1}{c}{\bf FMS} & \multicolumn{1}{c}{} & \multicolumn{1}{|c}{\bf (\%)} & \multicolumn{1}{c}{\bf FMS} & \multicolumn{1}{c|}{} & \multicolumn{1}{c}{\bf (\%)} & \multicolumn{1}{c}{\bf FMS} & \multicolumn{1}{c|}{} \\ \hline \hline
    \multicolumn{1}{|c|} {\multirow{4}{*}{\bf $\eta$ = 0.10}} &  \multicolumn{1}{c|}{\multirow{2}{*}{\bf $R$}}  & {\bf OPT} & 96.7 & 0.96  & 3.3e-1 & 100.0 & 1.00 & 9.1e-1 & 96.7 & 0.96 & 3.3e-1\\
    \multicolumn{1}{|c|} {} & \multicolumn{1}{c|}{} & {\bf ALS} & 100.0 & 0.99 &  & 100.0 & 1.00 &  & 90.0 & 0.90 & \\ \cline{2-12}
    \multicolumn{1}{|c|} {} & \multicolumn{1}{c|}{\multirow{2}{*}{\bf $ R+1$}}  & {\bf OPT} & 90.0 & 0.96 & 1.5e-13 & 100.0 & 1.00 & 6.8e-9 & 83.3 & 0.89 & 1.6e-11\\
    \multicolumn{1}{|c|} {} & \multicolumn{1}{c|}{} & {\bf ALS} & 3.3 & 0.24 &  & 13.3 & 0.52 &  & 10.0 & 0.13 & \\  \hline \hline
    \multicolumn{1}{|c|} {\multirow{4}{*}{\bf $\eta$ = 0.25}} & \multicolumn{1}{c|}{\multirow{2}{*}{\bf $R$}}  & {\bf OPT} & 76.7 & 0.97 & 3.3e-1  & 96.7 & 0.97 & 1.0 & 83.3 & 0.97 & 0.4e-1\\
    \multicolumn{1}{|c|} {} & \multicolumn{1}{c|}{} & {\bf ALS} & 73.3 & 0.94 &  & 96.7 & 0.97 &  & 70.0 & 0.84  & \\ \cline{2-12}
    \multicolumn{1}{|c|} {} & \multicolumn{1}{c|}{\multirow{2}{*}{\bf $R+1$}}  & {\bf OPT} & 86.7 & 0.97 &  6.4e-9 & 100.0 & 1.00 & 7.8e-5 & 76.7 & 0.90 & 1.4e-8 \\
    \multicolumn{1}{|c|} {} & \multicolumn{1}{c|}{} & {\bf ALS} & 13.3 & 0.40 &  & 46.7 & 0.72 &  & 10.0 & 0.23 & \\  \hline \hline
    \multicolumn{1}{|c|} {\multirow{4}{*}{\bf $\eta$ = 0.35}} & \multicolumn{1}{c|}{\multirow{2}{*}{\bf $R$}} & {\bf OPT} & 60.0 & 0.95 & 3.6e-1 & 100.0 & 1.00 & 3.3e-1 & 53.3 & 0.92 & 7.4e-1 \\
    \multicolumn{1}{|c|} {} & \multicolumn{1}{c|}{} & {\bf ALS} & 56.7 &  0.95 &  & 96.7 & 0.97 &  & 53.3 & 0.92 & \\ \cline{2-12}
    \multicolumn{1}{|c|} {} & \multicolumn{1}{c|}{\multirow{2}{*}{\bf $R+1$}}  & {\bf OPT} & 46.7 & 0.87 & 1.7e-9 & 100.0 &  1.00 & 4.0e-6 & 50.0 & 0.83 & 4.9e-7 \\
    \multicolumn{1}{|c|} {} & \multicolumn{1}{c|}{} & {\bf ALS} & 6.7 & 0.35 & & 40.0 & 0.62 &  & 10.0 & 0.30 & \\  \hline
\end{tabular}
\label{tab:ACC2}
\end{table*}


\section{Conclusions}
\label{sec:conclusions}
We have seen a shift in data mining in recent years: from models focusing on matrices to those studying higher-order tensors and now we are in need of models to explore and extract the underlying structures in data from multiple sources. One approach is to formulate this problem as a coupled matrix and tensor factorization problem. In this paper, we address the problem of solving coupled matrix and tensor factorizations when we have squared Euclidean distance as the loss function and introduce a first-order optimization algorithm called CMTF-OPT, which solves for all factor matrices in all data sets simultaneously. We have also extended our algorithm to data with missing entries by introducing CMTF-WOPT for the case where we have an incomplete higher-order tensor coupled with matrices. The algorithm can easily be extended to multiple incomplete data sets. 

To the best of our knowledge, the algorithms proposed so far for fitting a coupled matrix and tensor factorization model have all been alternating algorithms. We have compared our CMTF-OPT algorithm with a traditional alternating least squares approach and the numerical results show that all-at-once optimization is more robust to overfactoring as it has also been the case for fitting tensor models \cite{AcDuKo11a, ToBr06}.

Note that our current formulation of coupled matrix and tensor factorization has one drawback, which is encountered when we have a data set whose factor matrices are shared by all other data sets. For instance, we may have a tensor $\T{X}$, where $\T{X}=\KOp{\M{A}, \M{B},\M{C}}= \sum_{r=1}^R\lambda_r \MC{A}{r} \circ \MC{B}{r} \circ \MC{C}{r}$ and two matrices $\M{Y}=\M{A}\M{V}\Tra=\sum_{r=1}^R\alpha_r \MC{A}{r} \circ \MC{V}{r}$ and $\M{Z}=\M{B}\M{V}\Tra=\sum_{r=1}^R\beta_r \MC{B}{r} \circ \MC{V}{r}$. If we formulate a coupled factorization for these data sets as $ f(\M{A},\M{B},\M{C},\M{V})= \norm{\T{X} - \KOp{\M{A},\M{B},\M{C}}}^2 + \norm{\M {Y}-\M{A}\M{V} \Tra}^2 +\norm{\M {Z}-\M{B}\M{V} \Tra}^2$,
we do not take into account the cases for $\lambda_r \neq \alpha_r \neq \beta_r$. For such scenarios, scaling ambiguities should be taken into consideration by introducing a set of parameters for the scalars in the formulation.

Another issue we have not addressed in this paper is how to weigh different parts of the objective function modeling different data sets as discussed in \cite{Wilderjans2009}. This is an area of future research where a Bayesian framework for coupled matrix and tensor factorizations may be a promising approach. As another area of future research, we plan to extend our approach to different loss functions in order to be able to deal with different types of noise and incorporate nonnegativity constraints on the factors as nonnegativity often improves the interpretability of the model.

\section{Acknowledgments}
This work was supported in part by the Laboratory Directed Research and Development program at Sandia National Laboratories. Sandia National Laboratories is a multi-program laboratory managed and operated by Sandia Corporation, a wholly owned subsidiary of Lockheed Martin Corporation, for the U.S. Department of Energy's National Nuclear Security Administration under contract DE-AC04-94AL85000.

\bibliographystyle{abbrv}
\bibliography{paper}

\appendix

Here we briefly describe the data generation for Example 1 in \Sec{introduction}.
Tensor $\T{X}$ and matrix $\M{Y}$ coupled in the first mode are constructed as follows: 
\begin{compactitem}
\item \textbf{Step 1:} We form factor matrices
$\M{A}_1 \in \Real^{I \times R}$ and $\M{A}_2 \in \Real^{I \times R}$, where $R=2$, such
that in the first column of $\M{A}_1$, entries corresponding to the members
of $G_1$ and $G_2$ are assigned 1 (plus noise) while entries corresponding
to the members of $G_3$ and $G_4$ are assigned -1 (plus noise). It is
the vice versa for the second column; in other words, $G_3$ and $G_4$ members have 1 + noise values.
The columns of $\M{A}_2$ are generated similarly, except that in this case $G_1$ and $G_3$ form one
cluster while $G_2$ and $G_4$ form another.
\item  \textbf{Step 2:} Factor matrices $\M{B} \in \Real^{J \times R}, \M{C} \in \Real^{K \times R}$
and $\M{V} \in \Real^{M \times R}$ are generated using random entries following standard normal distribution.
\item \textbf{Step 3:} All columns of factor matrices are normalized to unit norm. $\T{X}$ and
$\M{Y}$ are constructed as $\T{X}=\KOp{\M{A}_1,\M{B},\M{C}}$ and $\M{Y}=\M{A}_2\M{V}\Tra$.
\end{compactitem}


\end{document}